\documentclass[10pt]{article}
\usepackage[dvips]{graphicx}
\usepackage{subfigure}
\usepackage{amsmath}
\usepackage{pstricks}
\usepackage{amssymb}
\usepackage{mathtools}
\usepackage{hyperref}
\usepackage{dsfont}

\newcommand{\sn}[0]{{\rm sn}}
\newcommand{\cn}[0]{{\rm cn}}
\newcommand{\dn}[0]{{\rm dn}}

\title{On the Alfred Gray's Elliptical Catenoid}
\author{Hugo Jim\'enez-P\'erez}
\date{\today}
\begin{document}
\maketitle

\begin{abstract}
    We give a parameterization, using Jacobi's elliptic functions, of 
    Alfred Gray's Elliptical Catenoid and Elliptical Hellicoid that
    avoids some problems present in the original depiction of these 
    surfaces.
\end{abstract}

\section{Introduction: The Bj\"orlings problem}
On of the Alfred Gray's favorite tools was the solution to Bj\"orling's
problem: \emph{Given a planar analytic curve, find a minimal surface in
$\mathbb R^3$ that contains it as a geodesic}. The Weierstrass representation 
provides an explicit solution: if $(x(t),y(t))$ is a parameterization of the
curve, the parameterization
\begin{eqnarray}
    \Phi(z) = \Re\left( x(z), y(z), i\int_p^z\sqrt{x^\prime(u)^2 +
        y^\prime(u)^2 }du \right)
    \label{eqn:bjor}
\end{eqnarray}
gives the solution, where $x(z)$, $y(z)$ are extensions of $x(t)$, $y(t)$ to
functions of a complex variable $z$. We will call it the Bj\"orling surface 
of the curve. When one takes, instead of the real part, the imaginary part of 
the expression, one gets another minimal surface, known as the conjugated
surface of the first one.

It is well know how Alfred used this to produce many beautiful surfaces. But
he could use it also as a theoretical tool: one day we asked Alfred what order
of contact  could a minimal surface have at a self-intersection point; without
a second of thought he replied: ``Any order of contact: take a planar curve
with a self-contact of order $n$ and solve Bj\"orling's problem''. An instant
theorem!.

\begin{figure}[h]
    \subfigure[]{
        \centering
	    \includegraphics[scale=0.5]{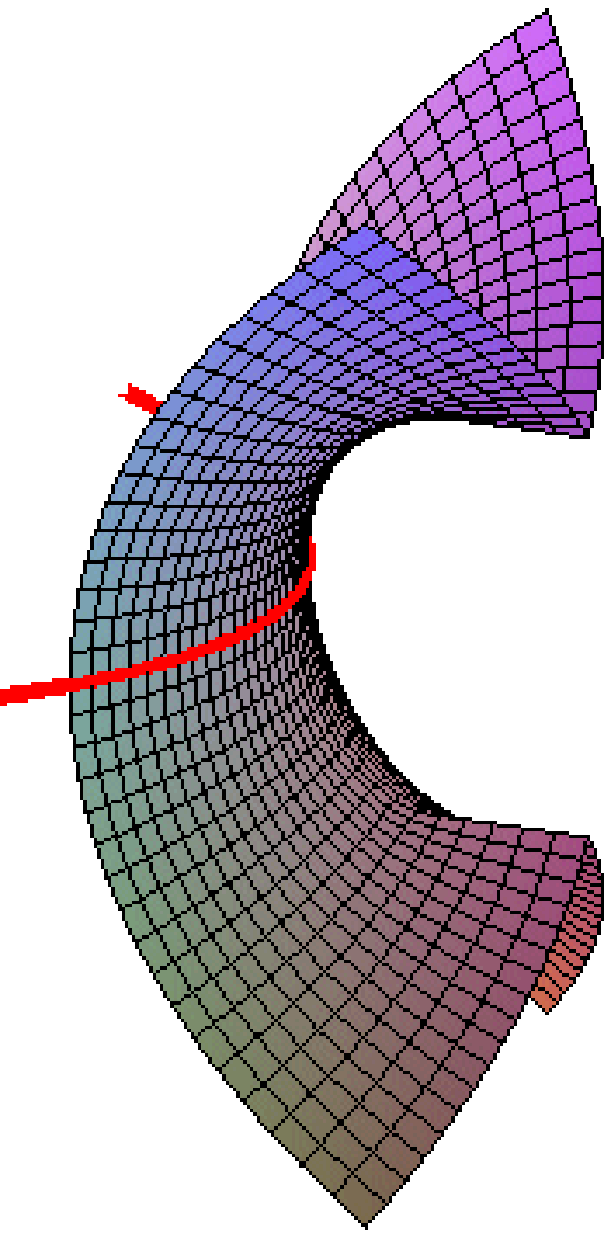}
    }
    \subfigure[]{
        \centering
	    \includegraphics{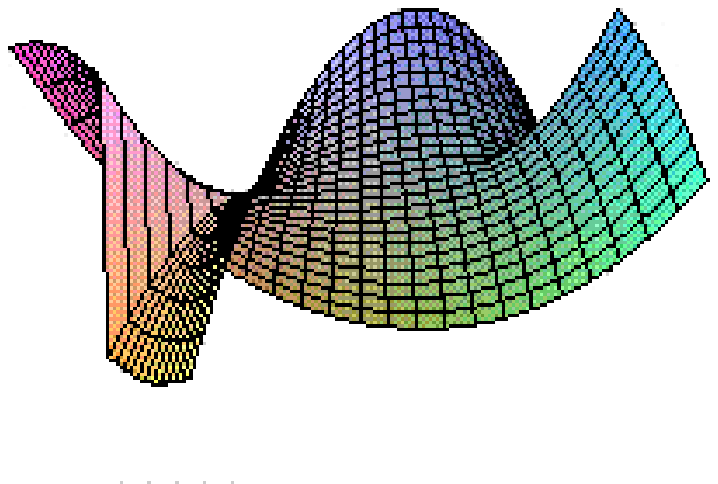}
    }
    \caption{The wrong picture for the parabola's Bj\"orling.}
    \label{fig:par_bad}
\end{figure}

Formula (\ref{eqn:bjor}) can be directly given to the computer. However, the
fact that the integrand is in many cases multivaluated poses some problems
when we try to get the global picture. Take, for instance, the parabola
$(2t, t^2)$. The formula gives the following parameterization of its Bj\"orling 
surface
\begin{eqnarray*}
    \Phi(z) &=& \Re\left( 2z, z^2, i\int_0^z\sqrt{1+u^2}du \right)
\end{eqnarray*}
which the reader can immediately give to his favorite graphics package\dots
and get the wrong picture! (Figure 1a, where the thick line shows the parabola).
The problem has to do with the integrand that is multivalued and branches 
at the points $i$ and $-i$. When we integrate from 0 to, say, $2i$ along
differents paths we can get different values, all with the same imaginary 
part. This produces a sharp edge, which is impossible in a minimal surface.
The problem is more evident in the conjugate surface: we get a discontinuity
which the computer fills with a planar face (Figure 1b).

\begin{figure}[h]
    \subfigure[Two views of Catalan's surface.]{
        \centering
	    \includegraphics[scale=0.7]{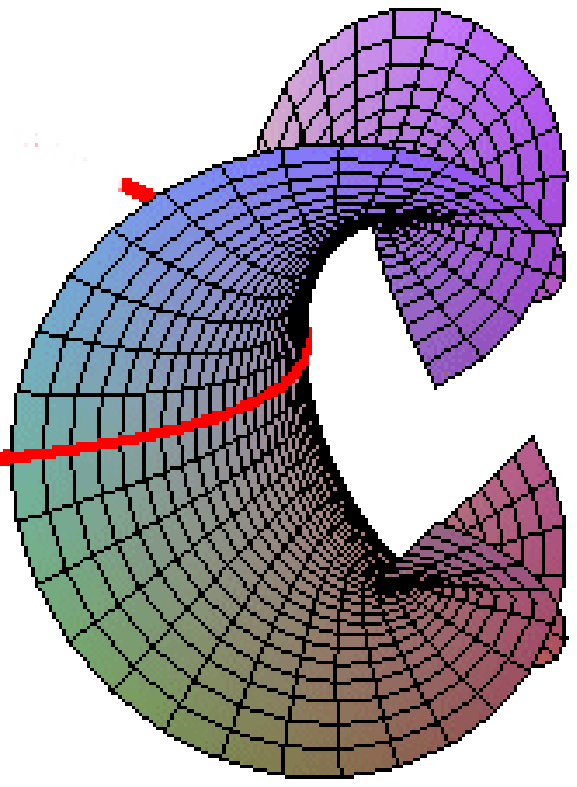}
	    \includegraphics[scale=0.5]{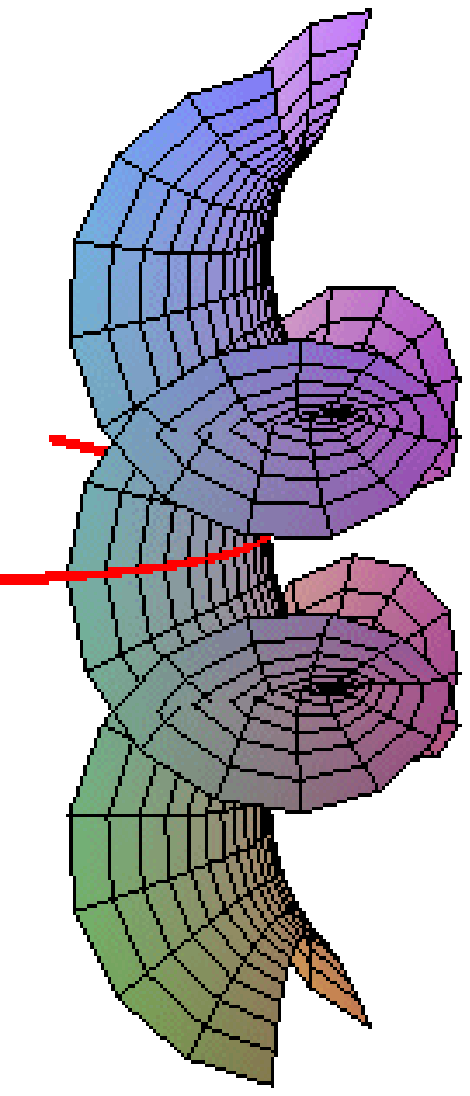}
    }
    \subfigure[The conjugate surface.]{
        \centering
	    \includegraphics[scale=0.5]{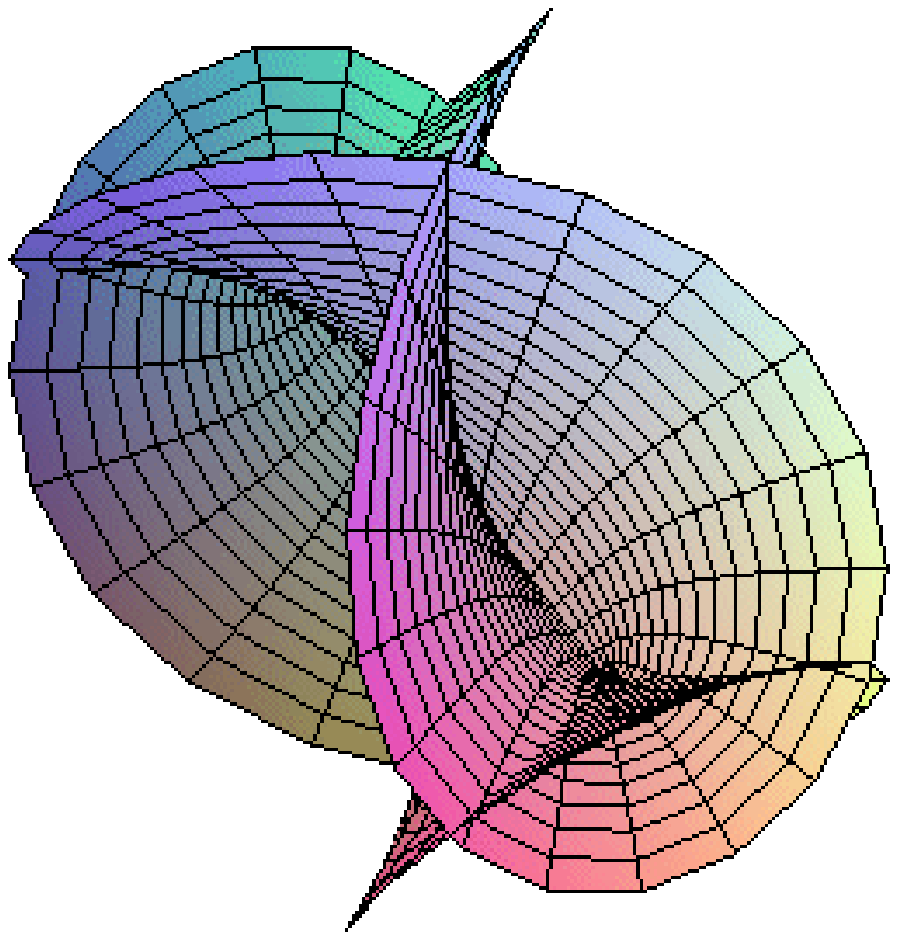}
    }
    \caption{The right picture for the parabola's Bj\"orling.}
    \label{fig:par_right}
\end{figure}

The problem can be easily solved by a simple sustitution $z=sinh w$ which
makes the function single valued and the surface becomes:
\begin{eqnarray*}
    \Phi(w) =\Re\left( 2\sinh(w), \sinh^2(w),2i\cosh(w\sinh(w)) + w \right)
\end{eqnarray*}
which gives the correct picture (Figure 2). The surface continues beyond  the
sharp edges which becomes lines of self-intersection and then takes a turn,
giving a periodic pattern. The surface turns out to be the well-know 
Catalan's surface that contains a cycloid perpendicular to the parabola
which is another planar geodesic of the surface. Figure 2b shows the correct 
conjugate surface.

\section{Bj\"orling's duality}
The above relation shows an interesting relation between the parabola 
and the cycloid: they have the same Bj\"orling surface. Recall that a 
plane cuts a minimal surface along a geodesic, if and only if, it is a plane of symmetry of the surface. Now, the Bj\"orling surface of a curve that has 
a line of symmetry has to planes of symmetry: the original one containing the 
curve and the plane perpendicular to it that contains the line of symmetry. In
formula (?) the original curve is obtained by restricting to real values 
of $z$ and the orthogonal geodesic by restricting to purely imaginary values
of $z$. The two curves are in a sense, dual to each other. To be more precise:
The duality is between objects consisting each of a curve and a point of
intersection of the curve with a line of symmetry. We can call those pairs
\emph{Bj\"orling duals} to each other. For example, the circle and the 
catenary are Bj\"orling duals (the common surface is the catenoid) and so
are the parabola and the cycloid with its line of symmetry that cuts
t on one of the highest (smooth) points of the cycle.

\section{The Elliptical Catenoid}
Let us apply (1) to the ellipse with semi-axes $a$ and $b$ $a>b$, and 
excentricity $e=\frac{\sqrt{a^2-b^2}}{b}$.
We  obtain 
\begin{eqnarray}
    \Phi(z) &=& \Re\left( b\cos(z), a\sin(z), i\int_0^z \sqrt{1-e^2\sin^2(u)}du \right)
    \label{eqn:cat1}
\end{eqnarray}
and with it the computer produces Figure 3. One observes that this surface has 
various \emph{channels} where one of the principal curvatures is close to zero
and the other one is very big. Again, this cannot be a minimal surface!
(One obtains a channel or an edge depending on the parity of the size 
of the grid used to plot the surface). The conjugate surface presents 
again a discontinuity which if filled by the computer with a flat face 
(Figure 3b)

Alfred's depiction in \cite{Gra1} of these surfaces, called by him 
\emph{Elliptical Catenoid} and \emph{Elliptical Hellicoid}, have the
same problems. (Those two loves of his, the computer and the Bj\"orling
problem, did not get togheter as smoothly as he thought). Both Maple
(which we use) and Mathematica give the same picture. Maybe Alfred was
aware of the problem since in this case he used an elliptic function 
instead of applying directly the formula. The analysis must, however, be
carried some steps further:

\begin{figure}[h]
    \subfigure[]{
        \centering
	    \includegraphics[scale=0.4]{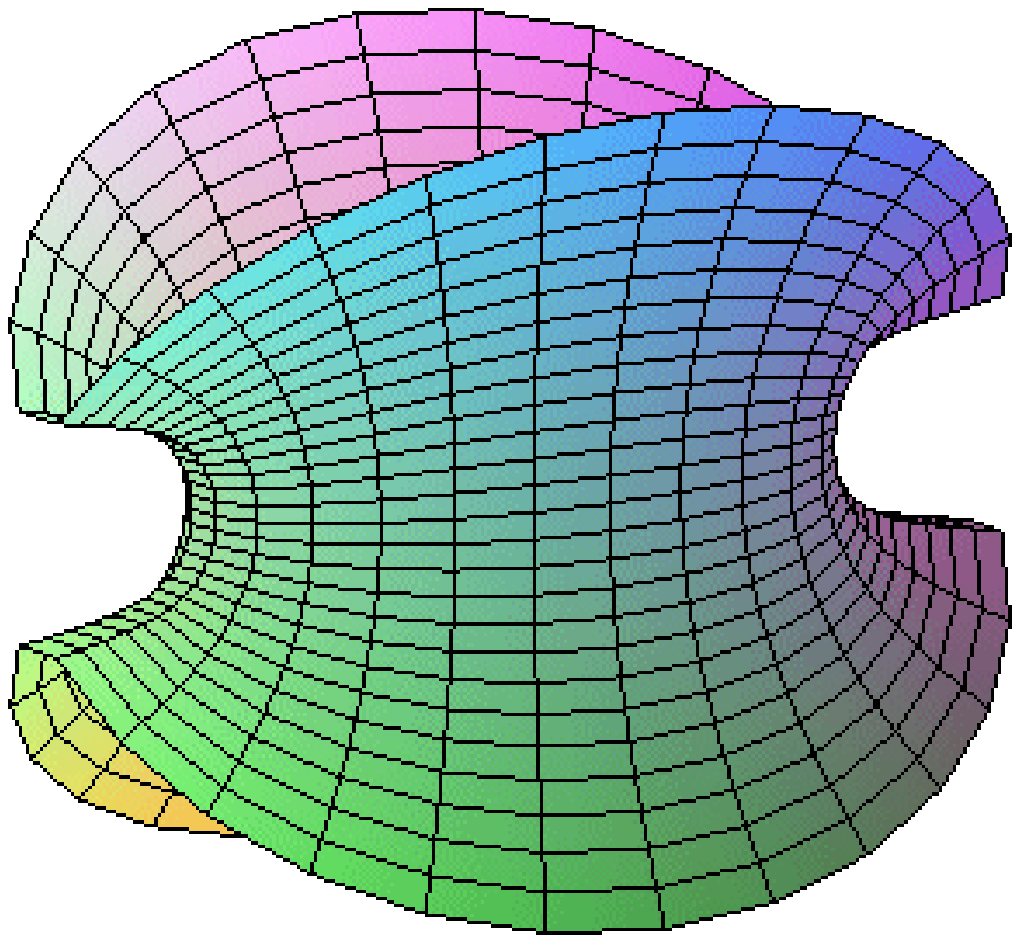}
    }
    \subfigure[]{
        \centering
	    \includegraphics[scale=0.8]{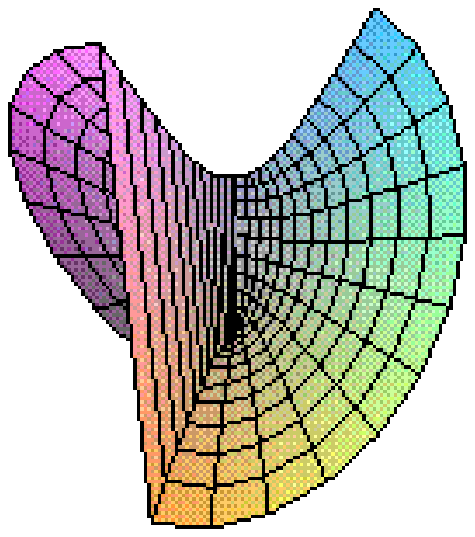}
    }
    \caption{The bad picture for the Elliptical Catenoid.}
    \label{fig:eli_bad}
\end{figure}

We use  Jacobi's elliptic functions (for details, see \cite{Law1}, \cite{Mar1}):
Jacobi's elliptic sine function $\phi=\sn(v,k)$, being the inverse function of
the integral 
\begin{eqnarray*}
    v = \int_0^\phi \frac{dt}{\sqrt{1-t^2}\sqrt{1-k^2t^2}},
\end{eqnarray*}
Jacobi's elliptic cosine
\begin{eqnarray*}
    \cn(v,k) = \sqrt{1-\sn^2(v,k)}
\end{eqnarray*}
and Jacobi's \emph{delta} function 
\begin{eqnarray*}
    \dn(v,k) = \sqrt{1-k^2\sn^2(v,k)}.
\end{eqnarray*}
Like the trigonometric funtions (which we obtain as the particular case $k=0$)
they can be extended to complex functions, but the extensions are now 
meromorphic and doubly periodic instead of being holomorphic and periodic.

With the new parameter $u$ related to $z$ by $\sn(u,e)=\sin(z)$ we obtain
the parameterization:
\begin{eqnarray*}
    \ Psi(u)=\Re\left( b\cn(u,e), a \sn(u,e),ia\int_0^u \dn^2(\sigma,e)d\sigma \right).
\end{eqnarray*}

The integrand is still a meromirphic function, but its residues are all $0$
\cite[p.241]{Law1}, so the integral is a well-defined function of $u$. This
parameterization produces the correct picture of the Elliptical Catenoid. 
The ellipse has two lines of symmetry and, correspondly, two Bj\"orling duals:
one resembles a cycloid, the other one a catenary. Figure 4 shows the correct 
pictures of the Elliptical Catenoid and the Elliptical Hellicoid:

\begin{figure}[h]
    \subfigure[The Elliptical Catenoid.]{
        \centering
	    \includegraphics[scale=0.4]{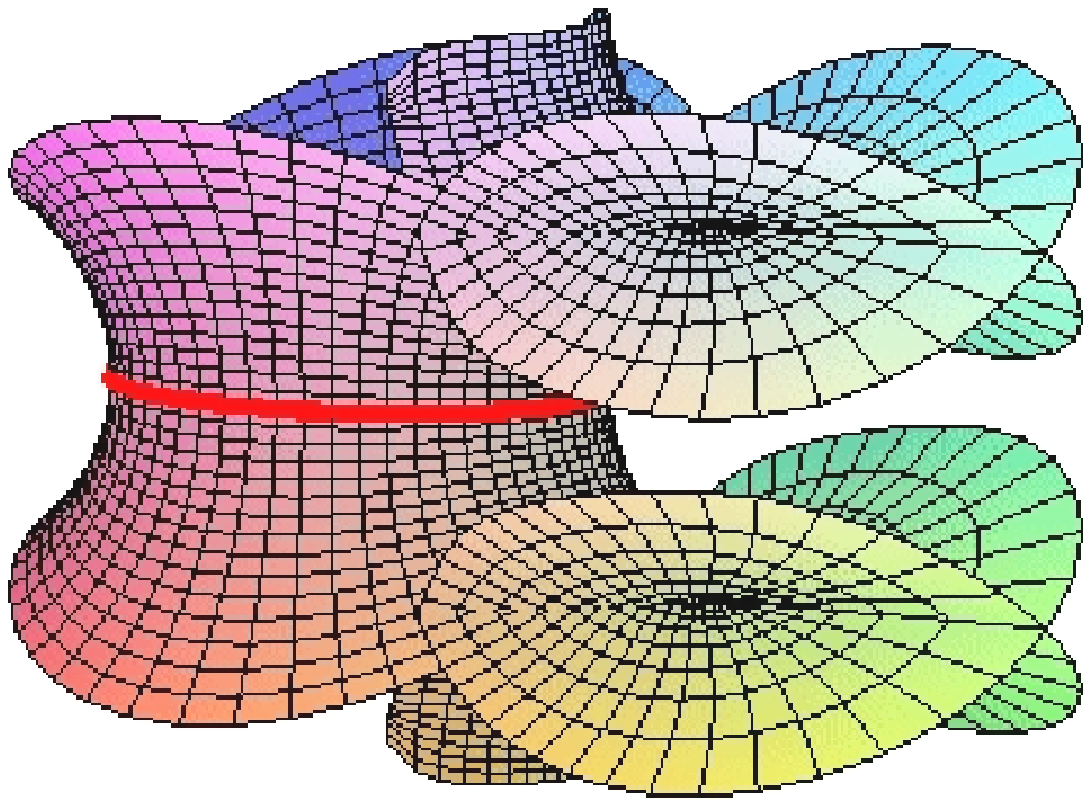}
    }
    \subfigure[The conjugate surface.]{
        \centering
	    \includegraphics[scale=0.4]{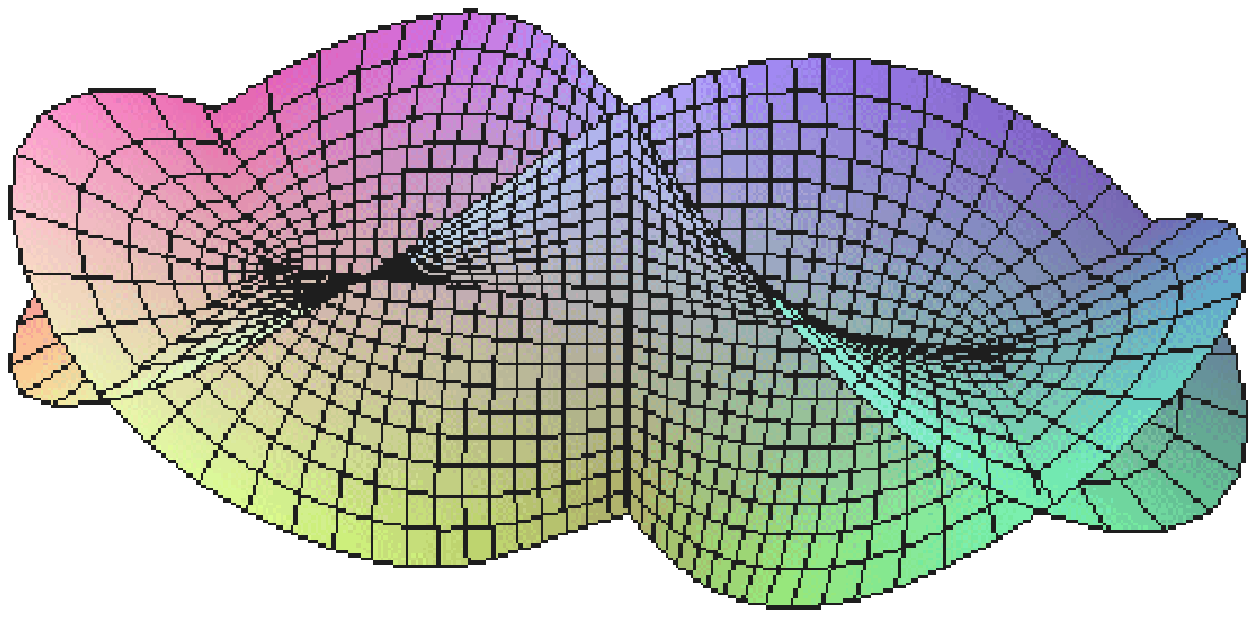}
    }
    \caption{The right picture for the Elliptical Catenoid.}
    \label{fig:eli_right}
\end{figure}
The same method gives the Bj\"orling surface of the hyperbola (where now we 
take $k=\frac{1}{e}<1$):
\begin{eqnarray*}
    \ Psi(u)=\Re\left( b\dn(u/k,k), iak \sn(u/k,k),a\int_0^u \cn^2(\sigma/k,k)d\sigma \right).
\end{eqnarray*}
The dual curve of the hyperbola again resembles a cycloid.

The analysis of the bj\"orling surfaces of the conics will be developed
fully in a forthcoming article.

\begin{figure}[h]
    \centering
    \subfigure[]{
        \centering
	    \includegraphics[scale=0.42]{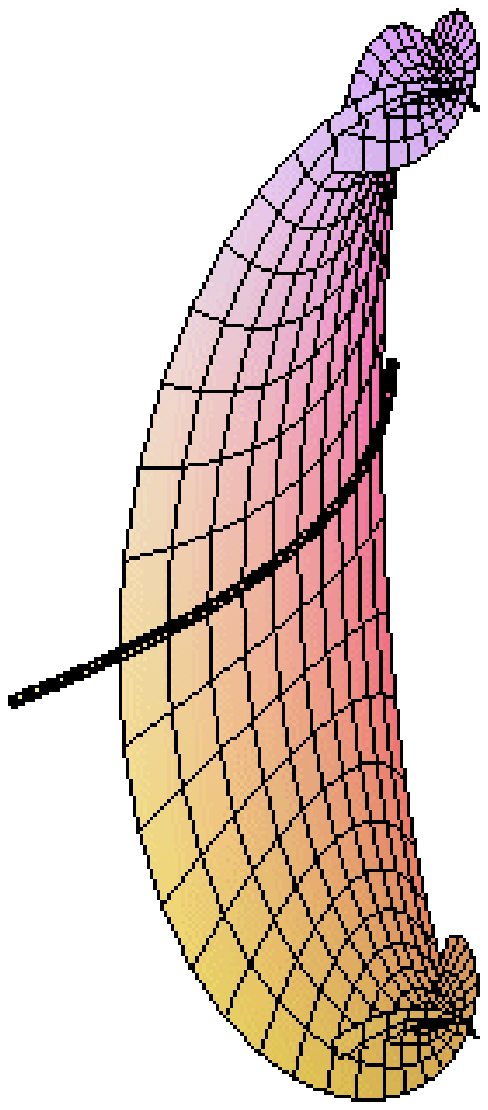}
	    }\hspace{80pt}
    \subfigure[]{
        \centering
	    \includegraphics[scale=0.5]{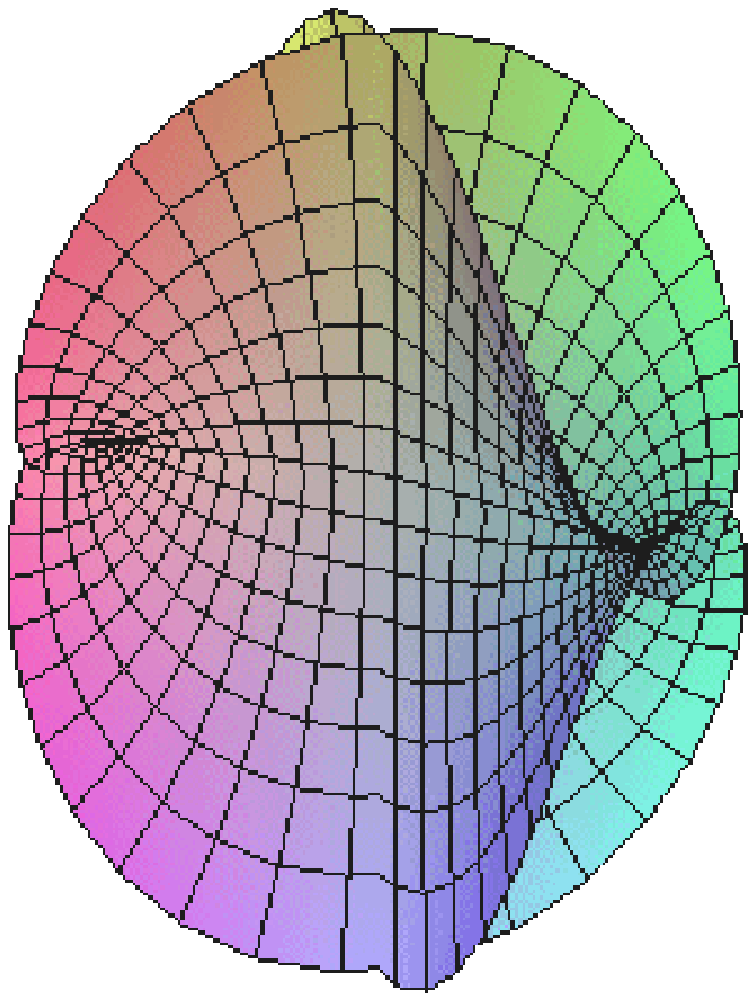}
    }
    \caption{The Bj\"orling surface for the hyperbola and its conjugated.}
    \label{fig:hyper_right}
\end{figure}

\end{document}